\documentstyle[12pt]{article}
\newcommand{\const}{\mathop{\rm const}\limits}

\newcommand{\mod}{\mathop{\rm mod}\limits}

\begin{document}

\begin{center}

{\bf A COUNTEREXAMPLE TO A HYPOTHESIS OF LIGHT TAIL }\\

\vspace{3mm}
{\bf OF MAXIMUM DISTRIBUTION FOR CONTINUOUS RANDOM } \\

\vspace{3mm}

 {\bf PROCESSES WITH LIGHT FINITE-DIMENSIONAL  TAILS.}\\

\vspace{4mm}

                E.Ostrovsky, \ L.Sirota.  \\

\vspace{2mm}
Department of  Mathematic, Bar-Ilan University,  Ramat Gan,
52900, Israel, \\
 e-mails: galo@list.ru,  \ eugostrovsky@list.ru ; \ sirota@zahav.net.il \\

\vspace{3mm}
                    {\sc Abstract.}\\

 \end{center}

 \vspace{3mm}

 We construct an example of a continuous centered random process with light tails of finite-dimensional  distribution
but with heavy tail of maximum distribution.\\

 \vspace{3mm}


\vspace{3mm}


\vspace{3mm}

\section{Notations. Statement of problem.}

 In the article \cite{Ostrovsky1}, 2008 year, has been formulated the following hypothesis $ H: $ \\

"Let  $ \theta = \theta(t), t \in T  $  be arbitrary separable
random field, centered: $  {\bf E} \theta(t) = 0, $  bounded with probability one:
$ \sup_{t \in T} |\theta(t)| < \infty \ (\mod {\bf P }), $ moreover, may be continuous, if  the set $ T $
is   compact metric space relative some distance. \par
Assume in addition that for some Young (or Young-Orlicz) function $  \Phi(\cdot) $ and correspondent Orlicz norm $ ||\cdot||Or(\Phi) $

$$
\sup_{t \in T}||\theta(t)||Or(\Phi) < \infty. \eqno(1)
$$

 Recall that the Luxemburg norm $ ||\xi||Or(\Phi)  $ of a r.v. (measurable function) $  \xi $  is defined as follows:

 $$
  ||\xi||Or(\Phi) = \inf_{k, k > 0} \left\{ \int_{\Omega} \Phi(|\xi(\omega)|/k) \ {\bf P}(d \omega) \le 1 \right\}.
 $$

  The Young function $ \Phi(\cdot) $ is by definition arbitrary even convex continuous strictly increasing on the
non-negative right-hand  semi-axis  such that
$$
\Phi(0) = 0, \ \lim_{u \to \infty} \Phi(u) = \infty.
$$

 Let also $ \Psi(\cdot) $ be {\it arbitrary} another Young function such that $ \lim_{u \to \infty} \Psi(u) = \infty, \ \Psi <<  \Phi,  $
which denotes by definition

$$
\forall \lambda > 0  \ \Rightarrow \lim_{u \to \infty} \frac{\Psi(\lambda u)}{\Phi(u)} = 0,  \eqno(2)
$$
see \cite{Rao1}, p.16. \par
 Recall that $ \Psi << \Phi $ implies in particular that the (unit) ball in the space $ Or(\Psi) $ is precompact set
in the space $ Or(\Phi). $ \par
\vspace{4mm}
{\it Open question: there holds (or not)}"

$$
|| \sup_{t \in T} \ |\theta(t)| \ ||Or(\Psi) < \infty. \eqno(3)
$$
\vspace{4mm}
 The conclusion (3) is true for the centered (separable) Gaussian fields \cite{Fernique1},
if the field $ \theta(\cdot) $
satisfies the so-called entropy or generic chaining condition \cite{Ostrovsky2},
\cite{Ostrovsky3},  \cite{Ostrovsky2}, \cite{Talagrand1}, \cite{Talagrand2}, \cite{Talagrand3},
\cite{Talagrand4}; in the case when  $ \theta(\cdot) $ belongs to the domain of attraction of
Law of Iterated Logarithm \cite{Ostrovsky5} etc.\par

 Notice that if the field $ \theta(t) $ is continuous $ ( \mod {\bf P}) $  and satisfies the condition
(1), then {\it there exists } an Young function $ \Psi(\cdot), \ \Psi(\cdot) << \Phi(\cdot) $
for which the inequality  holds, see \cite{Ostrovsky1}. \par

The condition of a view $ ||\xi||Or(\Psi) < \infty $ described the tail behavior for
the distribution of the random  variable $ \xi. $  Another approach which was used in the
monograph M.Ledoux and M.Talagrand \cite{Talagrand1}, p. 309-317 is related in fact with generalized
Lorentz  (more exactly, Lorentz -Zygmund) norm $ ||  \ \xi \ ||L(v):  $

$$
 || \xi ||L(v) \stackrel{def}{=} \sup_{ A: {\bf P}(A)> 0}  \left[ \frac{1}{v({\bf P}(A))} \cdot \int_A |\xi(\omega)| \ {\bf P}( d \omega) \right].
$$
 Here $ v = v(z), \ z \in (0,1] $ is continuous monotonically increasing function   such that  $ v(0) = 0. $ \par

 Notice that in all this cases  the inequality (3) is true with replacing the function $  \Psi $
on the function $ \Phi.$\par

\vspace{3mm}
{\bf  Our target is to give a negative answer on this  hypothesis, by means of construction of
correspondent counterexample. }\\

\vspace{2mm}

{\bf  We consider in this short article also the case  of unbounded measure   } {\bf P.  } \\

\vspace{3mm}

  The detail investigation of the theory of Orlicz's spaces included the case unbounded measure see
in the monographs \cite{Rao1}, \cite{Rao2}.\par

 The notion of separability and  entropy or generic chaining condition  for the conclusion (3)
in this case see in  \cite{Ostrovsky7}. \par

\vspace{4mm}

 {\sc  Several notations and definitions.}\par

\vspace{3mm}

 {\bf A.} A triplet $ (\Omega, \cal{B}, {\bf P} ),  $  where $ \Omega = \{\omega\} = \{x\}  $ is arbitrary set,
$ \cal{B} $ is non-trivial sigma-algebra  subsets $ \Omega  $ and  $ {\bf P}  $ is non-zero non-negative completely
additive measure defined on the $  \cal{B} $ is called a probabilistic space, even in the case when $ {\bf P}(\Omega) = \infty. $
   We denote as usually for the random variable $ \xi $ (r.v.) (i.e. measurable function $ \Omega \to R )  $

 $$
 |\xi|_p = [ {\bf E} |\xi|^p ]^{1/p} = \left[ \int_{ \Omega  } |\xi(\omega)|^p \ {\bf P}(d \omega)  \right]^{1/p}, \ p \ge 1;
 $$

$$
L_p = \{ \xi, \ |\xi|_p < \infty  \}.
$$

\vspace{3mm}

{\bf B.} The so-called Grand Lebesgue Space $ G \psi $  with norm $  ||\cdot ||G\psi   $ is defined (in this article) as follows:

$$
G \psi = \{\xi, \ ||\xi||G \psi < \infty \}, \  ||\xi||G \psi \stackrel{def}{=} \sup_{ p \ge 1 } \left[ \frac{|\xi|_p}{\psi(p)} \right].
\eqno(4)
$$
 Here $ \psi = \psi(p) $ is some continuous increasing function such that $ \lim_{p \to \infty} \psi(p) = \infty. $\par
 The detail investigation of this spaces (and more general spaces) see in \cite{Liflyand1},  \cite{Ostrovsky8}. See also
 \cite{Fiorenza1}, \cite{Fiorenza2}, \cite{Iwaniec1}, \cite{Iwaniec2}, \cite{Kozachenko1} etc. \par
\vspace{4mm}
   {\it An important for us fact about considered here spaces is proved in \cite{Ostrovsky6}: they coincide with some exponential  Orlicz's spaces}
  $ Or(\Phi_{\psi}). $ For instance, if $ {\bf P}(\Omega) = 1 $ and  $ \psi(p) = \psi_{1/2}(p): =\sqrt{p}, $ then the space
  $ G\psi_{1/2} $ consists on all the subgaussian (non-centered, in general case) r.v. $ Or(\Phi_{\psi_{1/2}})  $
   for which $ \Phi_{\psi_{1/2}}(u) = \exp(u^2/2) - 1.  $  \par
 The Gaussian distributed r.v. $ \eta $ belongs to this space.  Another example: let $ \Omega = (0,1) $ with usually Lebesgue measure and

 $$
 f_{1/2}(\omega) = \sqrt{|\log \omega}|, \ \omega > 0; \ f_{1/2}(0)= 0.
 $$
  It is easy to calculate using Stirling's formula for the Gamma function:
 $$
 | f_{1/2}|_p \asymp \sqrt{p}, \ p \in (1,\infty).
 $$
  The tail behavior:

  $$
  {\bf P} ( f_{1/2} > u ) = \exp (-u^2).
  $$

 The case when in (4) supremum is calculated over {\it finite } interval is investigated in
 \cite{Liflyand1}, \cite{Ostrovsky8}, \cite{Ostrovsky9}:

$$
G_b \psi = \{\xi, \ ||\xi||G_b \psi < \infty \}, \  ||\xi||G_b \psi \stackrel{def}{=} \sup_{ 1 \le p < b } \left[ \frac{|\xi|_p}{\psi(p)} \right],
 \  b = \const > 1,  \eqno(5)
$$
but in (5) $ \psi = \psi(p) $ is  continuous  function in the semi-open interval $ 1 \le p < b $
 such that $ \lim_{p \downarrow b} \psi(p) = \infty. $\par

 An used further example:

 $$
 \psi^{(\beta,b)}(p) = (b-p)^{-\beta}, \ 1 \le p < b, \beta = \const > 0; \ G_{\beta,b}(p) := G_b\psi^{(\beta,b)}(p).
 $$

\vspace{3mm}

{\bf C.} Recall that sets $ A_1, A_2, \ A_i \in \cal{B} $ are disjoint,  if $ A_1\cap A_2 = \emptyset. $ The sequence of a functions
$  \{h_n \}, n =1,2,3  \ldots  $ is said to be {\it disjoint}, or more exactly {\it pairwise disjoint, if }

$$
\forall i,j; i \ne j \ \Rightarrow \ h_i \cdot h_j = 0. \eqno(6)
$$
If the sequence of a functions $  \{h_n \} $  is pairwise disjoint, then

$$
|\sum_n h_n|_p^p = \sum_n |h_n|_p^p, \hspace{6mm}  \sup_n |h_n(x)| = \sum_n |h_n(x)|. \eqno(7)
$$

\vspace{3mm}

{\bf D.} We denote as ordinary  for any measurable set $ A, \ A \in \cal{B} $ it indicator function by
$ I(A) = I_A(\omega).   $\par

\vspace{3mm}

\section{Main result.}

\vspace{3mm}

{\bf Theorem.}  The proposition of hypothesis $ H $  is not true. Indeed. there exist: \\

{\bf A.} A centered  continuous in the $ Or(\Phi) $ sense and with probability one
 random process (field) $ \theta = \theta(t) = \theta(t,\omega) $
defined aside from the probabilistic space on some {\it compact } metric space $ (T,d) = (\{t \}, d), $
where $ d $ is the non-trivial distance in the set $ T, $
such that for some Young function $ \Phi = \Phi(u) $ the condition (1) is satisfied.  \par

\vspace{4mm}

{\bf B.} An Young function $ \Psi = \Psi(u) $ for which $  \lim_{u \to \infty} \Psi(u) = \infty, \ \Psi << \Phi $   but

\vspace{4mm}

{\bf C.}
$$
  || \sup_{t \in T} \ |\theta(t)| \ ||Or(\Psi) = \infty. \eqno(8)
$$

{\bf Proof.}\\

\vspace{3mm}

{\bf 1.} In the sequel  we choose as a capacity metric space $  (T,d) $ the set of positive integer
numbers with infinite point which we denote $  \infty: $

$$
T = \{ 1,2,3, \ldots, \infty   \}.
$$
 The distance $  d $ is defined as follows:
 $$
 d(i,j) = \left| \frac{1}{i} - \frac{1}{j} \right|, \ i,j < \infty; \ d(i,\infty) = d(\infty,i)= \frac{1}{i}, \ i < \infty; \eqno(9)
 $$
and obviously $ d(\infty,\infty) = 0.  $\par
 The pair $ (T,d) $ is compact (closed) metric space  and the set $ T $ has unique limit point $ t_{0} = \infty. $
For instance, $ \lim_{n \to \infty} d(n,\infty) = 0.   $ \par

\vspace{4mm}

{\bf 2.}  We consider  first of all the case $ {\bf P} (\Omega) = \infty.  $  Introduce as an example
the following triplet: $ \Omega = \{ x \} = R_+, \ \cal{B}  $ is Borelian sigma-algebra and $ {\bf P} $
is Lebesgue measure: $  {\bf P} (dx) = dx. $ \par
 Let us define a numerical sequences $ c(n) = \log^{-3}(n + 3), \ n \in T $ and a sequence of a functions

$$
g_n(x) = c(n) \cdot I_{(n, n+1)}(x) \cdot f_{1/2}(x-n), \ g_{\infty}(x) = 0,  \eqno(10)
$$

$$
g(x)  = \sum_{n=1}^{\infty} g_n(x). \eqno(11)
$$
Note that the functions  $ g_n $ are disjoint and following $ \sup_n |g_n(x)| < \infty  $ almost surely. \par
We calculate using the relations (7):

$$
|g_n|_p = c(n) \ \psi_{1/2}  (p),  \eqno(12)
$$
therefore  $  ||g_n - g_{\infty}||G\psi_{1/2} = ||g_n ||G\psi_{1/2} \to 0  $ as $ n \to \infty $ and moreover
$ g_n \to 0  $ almost everywhere.  Indeed, let $ \epsilon $ be arbitrary positive number. We get:

$$
\sum_n {\bf P} ( |g_n| > \epsilon ) \le \sum_n \exp (-C \epsilon^2 \log^6(n) ) < \infty.
$$
 Our conclusion follows from the lemma of Borel- Cantelli, which is true even for unbounded measure. \par

  So, the process $ \theta(t) = g_n, \ n = t  $ satisfies the condition (1)
relative the Young function $ \Phi_{\psi_{1/2}}(u). $\par
 Let now $ p $ be arbitrary number,  $ 1 \le p < \infty. $  We have:

 $$
  | \ \sup_n |g_n| \ |_p^p =  \sum_n |g_n|_p^p    = \sum_n c^p(n) \psi^p_{1/2}(p) = \infty.\eqno(13)
 $$

  So, we can choose in the capacity  of the Young  function $ \Psi = \Psi(u) $ the function  $  \Psi(u) = |u|^p. $ \par

  In order to obtain the centered needed process $ \theta(t)  $  with at the same properties,
 we consider the sequence  $  \tilde{g}_n(x) = \epsilon(n) \cdot g_n(x),  $  where $ \{\epsilon(n) \} $ is a Rademacher
 sequence independent on the $  \{g_n\}: $

$$
{\bf P}(\epsilon(n) = 1) = {\bf P}(\epsilon(n) = -1) = 1/2;
$$
then  $  |\tilde{g}_n(x)| = |g_n(x)|,  \ | |\tilde{g}_n|_p = |g_n|_p $
 and the sequence $ \{\tilde{g}_n \} $ is also pairwise disjoint (Rademacher's symmetrization).\par
  This completes the proof of our theorem, but only in the case $ {\bf P} (\Omega) = \infty. $\par

\vspace{4mm}

{\bf 3.} The case $  {\bf P} (\Omega) = 1 $  is more complicated. We choose $ \Omega = (0,1)  $ with
Lebesgue measure  and define

$$
a(n) = 1 - 0.5  n^{-\alpha}, \  \alpha = \const \in (0,1); \ p_0 = \const > 1;
$$

$$
\Delta(n) = a(n+1) - a(n) \sim_{n \to \infty} C_{\alpha} \ n^{-\alpha -1};  \ c(n) = n^{\alpha/p_0}.  \eqno(14)
$$
 We introduce as before the following  positive random process $ \theta(t) = g_n, \ n = t, \ t,n \in T, \ \Omega = \{x \},  $

$$
g_n(x) = c(n) \ f_{1/2} \left( \frac{x-a(n)}{\Delta(n)} \right) \ I_{(a(n), a(n+1)) }(x), \ x \in \Omega, \
g_{\infty}(x) = 0;   \eqno(15)
$$

$$
g(x) = \sum_n g_n(x) =  \sum_n n^{\alpha/p_0 } \ f_{1/2} \left( \frac{x-a(n)}{\Delta(n)} \right) \ I_{(a(n), a(n+1)) }(x).  \eqno(16)
$$

 Note that the sequence of r.v.  $ \{ g_n(x) \}  $ is again non-negative and  disjoint, therefore

 $$
 \sup_n g_n(x) = \sum_n g_n(x) = g(x), \ |\sup_n g_n|_p^p = \sum_n |g_n|_p^p. \eqno(17)
  $$

\vspace{3mm}

 {\bf Remark 1.} In this pilcrow  $ c(n) \to \infty, $ in contradiction to the case $ {\bf P } (\Omega) = \infty. $\par

\vspace{3mm}

  Let now $ p \in [1, p_0);$  in what follows we presume $ p \to p_0-0 .$   We calculate consequently,
 as long as the function $  \psi_{1/2}(p) $ is bounded: $ \psi_{1/2}(p) \le \sqrt{p_0}: $

$$
|g_n|_p^p = n^{\alpha p/p_0} \ \Delta(n) \sim C_{\alpha} n^{-1 - \alpha(1-p/p_0)},  \eqno(18)
$$
therefore $ |g_n - g_{\infty}|_{p_0} \to  0, \  n \to \infty. $ The $  L_{p_0} $ and with probability one
continuity of the process $ g_n, \ n \in T $
it follows from  (18) and the main result  of an article of G.Pizier \cite{Pizier1}. \par
 Further,

$$
|g|_p^p/C_{\alpha} = \sum_n  n^{-1 - \alpha(1-p/p_0)}   \sim C_1(\alpha) \ (p_0 -p)^{-1};
$$

$$
|g|_p \sim C_2(\alpha) \ (p_0 -p)^{ -1/p } \sim C_3(\alpha) \ (p_0 -p)^{ -1/p_0 }. \eqno(19)
$$

 The equality (19) implies on the language of $ G\psi $ spaces that

 $$
|| \ \sup_n \ |g_n| \ ||G_{1/p_0,p_0} < \infty,  \eqno(20)
 $$
and that the relation (20) is exact. \par

 It follows from the equalities (19) and (20) that the tail function for the r.v. $ g = g(x), $
i.e. the function $ G_g(z) := {\bf P}(|g(x)| \ge z), z > 0  $  obeys the following asymptotical
 as $  z \to \infty $ expression:

$$
 G_g(z) \sim C_4(\alpha) \ z^{-p_0}, \eqno(21)
$$
i.e. at the same asymptotic as for the r.v. $ \eta =  \eta(x) = C_5(\alpha) \ x^{ -1/p_0 }. $\par
 We offer in the capacity of the Young function $ \Phi = \Phi(u) $ the Young-Orlicz function for the $L_{p_0}  $
space: $ \Phi(u) = |u|^{p_0}. $ \par
 Let us consider the following Young function:

 $$
 \Psi(u) = |u|^{p_0} \ ( \log ( e + |u| ))^{ -1/2 }. \eqno(21)
 $$
  By virtue of the condition $ {\bf P} (\Omega) = 1  $ we deduce $  \Psi << \Phi.$ \par

 Since the Orlicz spaces  are rearrangement invariant,

 $$
 ||g||Or(\Psi) < \infty \ \Leftrightarrow ||\eta||Or(\Psi) < \infty.
 $$

 As long as the function $ \Psi =\Psi(u) $ satisfies the $ \Delta_2 $ condition,
$$
||\eta||Or(\Psi) < \infty \Leftrightarrow  \int_{\Omega} \Psi(|\eta(x)| \ dx < \infty.
$$

 But

 $$
 \int_{\Omega} \Psi(|\eta(x)| \ dx  \ge C_6 \ + C_7 \ \int_0^{1/2}  x^{-1} \ |\log x|^{-1/2} \ dx = \infty.  \eqno(22)
 $$

 Thus,

 $$
 || \ \sup_n |g_n| \ || G\Psi  = || \ |g| \ || G\Psi = \infty. \eqno(23)
 $$

 It remains to use the known Rademacher's symmetrization method in order to obtain the centered  process $ \theta(t). $ \par
\vspace{3mm}

{\bf Remark 2.} The constructed process $ \theta(t) $ give us a new example of centered continuous random process with light
tails of finite-dimensional distribution, but for which  entropy and generic chains series divergent.\par

\vspace{3mm}
{\bf Remark 3.} The proposition of our theorem remains true with at the same  (counter) example
if we use instead Orlicz space the generalized Lorentz   (Lorentz -Zygmund)
space $ L(v) $ or the space  $ K = K(h) = \{\tau = \tau(\omega) \} $
with quasinorm
$$
||\tau ||K = \sup_{z > 0} [ G_{\tau}(z)/ h(z)].
$$
 Here $ h = h(z), \ z \in (0,\infty) $ is continuous monotonically decreasing function   such that  $ v(0+) = {\bf P} (\Omega). $ \par
It is important for us only that the tail of distribution for $ \sup_t \theta(t) $ is essentially greatest in comparison with
$ 1/\Phi(u/C) $  for arbitrary constant $ C > 0. $\par

\vspace{3mm}
{\bf Remark 4.} The proposition of our theorem remains true if we use instead the space of continuous function $ C(T,d)  $
arbitrary  separable Banach space.\par

\end{document}